\documentclass[11pt,reqno]{amsart}
\usepackage{amssymb}
\usepackage[all]{xy}

\newtheorem{thm}{Theorem}[section]
\newtheorem{lem}[thm]{Lemma}
\newtheorem{cor}[thm]{Corollary}
\newtheorem{prop}[thm]{Proposition}
\newtheorem{hyp}[thm]{Assumptions}

\theoremstyle{definition}
\newtheorem{rem}[thm]{Remark}

\newtheorem{exa}[thm]{Example}

\newfont{\cyrr}{wncyr10}

\def\x{\times}
\def\ox{\otimes}
\def\o+{\oplus}
\def\ra{\rightarrow}
\def\ua{\uparrow}
\def\lra{\longrightarrow}

\def\da{\downarrow}
\def\beqa{\begin{equation*}}
\def\eeqa{\end{equation*}}

\def\Ext{\rm Ext}

\def\al{\alpha}

\def\la{\lambda}
\def\si{\sigma}

\def\C{{\mathcal C}}

\def\F{\mathcal F}

\def\cO{{\mathcal O}}

\title{ STABLE BUNDLE EXTENSIONS ON ELLIPTIC CALABI-YAU THREEFOLDS}

\author{Bj\"orn Andreas}
\address{Institut f\"ur Mathematik und Informatik,
Freie Universit\"at Berlin, Arnimallee 14, 14195 Berlin, Germany}
\email{andreas\char`\@math.fu-berlin.de}

\author{Gottfried Curio}
\address{Arnold-Sommerfeld-Center for Theoretical Physics,
Department f\"ur Physik, Ludwig-Maximilians-Universit\"at M\"unchen,
Theresienstr. 37, 80333 M\"unchen, Germany}
\email{curio\char`\@theorie.physik.uni-muenchen.de}

\thanks{B. A. is supported by DFG-SFB 647/A1. Report-no: LMU-ASC 76/06.}

\begin{document}
\begin{abstract}
We construct stable bundle extensions on elliptically fibered Calabi-Yau
threefolds. We show that these bundles can solve the topological anomaly constraint
in heterotic string theory without the need of invoking background fivebranes. 

\end{abstract}

\maketitle

\section{Introduction}

In the present paper we consider the class of elliptically fibered Calabi-Yau threefolds $\pi\colon X\to B$ with a section $\si$ and construct stable vector bundles $V$ of vanishing first Chern class using the method of bundle extensions. For a choice of data the bundles satisfy 
\beqa
c_2(TX)=c_2(V)\nonumber
\eeqa
and so qualify as physical gauge bundles for heterotic string theory compactification.

In contrast, the widely used spectral cover construction \cite{FMW}, \cite{FMWIII}, \cite{Don1} gives stable 
vector bundles on elliptic fibrations. These bundles solve the generalized anomaly constraint $c_2(TX)-c_2(V)=[W]$ with $[W]$ an effective curve class
(cf. below, Section 7). This mismatch causes two problems in physical model building: first, for $[W]\neq 0$ it prevents 
the model to be interpreted as a non-linear sigma model; secondly, for $[W]\neq 0$, and even for $[W]=0$ while $V\neq TX$, it leads to singular 3-forms in the heterotic anomaly equation. As a consequence, it is more difficult to solve the anomaly equation because a non-trivial string theory $H$-field has to be taken into account.  A detailed discussion of these issues is given in the physical companion paper \cite{ACu06}.

The paper is organized as follows. In Section 2, a general outline of the construction method is given and the necessary steps for proving stability of a given non-split extension are described. We consider extensions of $W_q \ox\cO_X(pD)$ by $U_p\ox \cO_X(-qD)$ where $U_p$ and $W_q$ are given stable vector bundles of vanishing first Chern class. $D$ is a divisor in $X$, chosen such that the resulting vector bundle has trivial determinant. The main result of this section is Lemma \ref{21}, which gives a sufficient condition for the bundle not to be destabilized by certain subsheaves. In Section 3, we restrict to extensions $V$ of $\cO_X(nD)$ by $\pi^*E\ox \cO_X(-D)$ 
where $E$ is a given stable vector bundle on the base of the Calabi-Yau threefold. The main result of this section is Lemma \ref{nspE}; it provides conditions for choosing such an extension $V$ non-split. In Section 4, we prove the stability of a pull-back bundle $\pi^*E$, assuming the base of the Calabi-Yau space is given by the Enriques surface (Lemma \ref{pE}). We then prove that $V$ is stable in a specific region of the K\"ahler cone (Proposition \ref{thm22}). In Section 5, we generalize the results of the previous section to elliptic Calabi-Yau threefolds with ample $K_B^{-1}$. In Section 6, we consider extensions of stable spectral cover bundles $V_n$ and prove stability for extensions of  $\cO_X(n\pi^*\al)$ by $V_n\ox \cO_X(-\pi^*\al)$ where $\al\in H^2(B, {\mathbb Z})$. In Section 7, we give explicit solutions to the topological anomaly constraint imposed by heterotic string theory.

Throughout the paper we use the notation $c_i=c_i(B)$.
\\ \\ 
B. A. would like to thank D. Hern\'andez Ruip\'erez and H. Kurke for helpful discussions.


\section{Method of construction}

Let $\pi\colon X\to B$ be an elliptically fibered Calabi-Yau threefold with a section $\si$. Except for Section 6, we will consider $B$ to be either an Enriques surface or a surface with ample $K_B^{-1}$ such as 
 the Hirzebruch surface ${F_r}$ with $r=0,1$ or the del Pezzo surface $dP_k$ with $k=0,\dots, 8$. 
 
We consider vector bundles $V$ of rank $p+q$ on $X$ defined as non-trivial extensions of bundles $U_p$ and $W_q$.

\begin{hyp}
We assume that $U_p$ and $W_q$ are stable vector bundles both of vanishing first Chern class.
\end{hyp}

Furthermore, $U_p$ and $W_q$ should be suitably twisted by line bundles such that $V$ has vanishing first Chern class
$$
0 \to  U_p\ox \cO_X(-qD)\to  V \to W_q \ox\cO_X(pD) \to  0
$$
where $D=x\si + \pi^*\al$. To discuss stability we will choose as polarization $J=z\si + \pi^*H$
where $H$ (chosen in the integral cohomology)
is in the K\"ahler cone $\C_B$ of the base $B$ and $z\in {\mathbb R^{>0}}$.
For an elliptically fibered Calabi-Yau space $X$ one finds $J$ in the K\"ahler cone $\C_X$ of $X$ under the following conditions \cite{Bfibre}
\begin{equation}
J\in \C_X \Longleftrightarrow  z>0 \; , \;\; H - z c_1 \in \C_B.\nonumber
\end{equation}
The following Lemma is obvious.
\begin{lem}
\label{obv}
Two necessary conditions for $V$ to be stable are
\begin{enumerate}
\item $DJ^2>0$ or equivalently $\mu(U_p\ox \cO_X(-qD))<0$
\item $W_q \ox\cO_X(pD)$ is not a subbundle of $V$, i.e. the exact sequence defining $V$ can be chosen non-split.
\end{enumerate}
\end{lem}
For the rest of this section let us assume that the non-split condition can be satisfied.

To discuss the required steps for proving stability of $V$ consider the following diagram of exact sequences \[ \begin{array}{ccccccccc}
  &  &0 & & 0   &  &0 & \\
   &      & \ua &                   &\ua&             & \ua &      &\\
0 &\to &  P=\bar{P}\ox \cO_X(-qD) & \to & V/V'_{r+s}&\to & T=\bar{T}\ox\cO_X(pD)&\to &0\\
  &      & \ua &                   &\ua&                    & \ua &      &\\
0 & \to & U_p\ox \cO_X(-qD)& \stackrel{i}{\to} & V & \stackrel{j}{\to}
 &W_q \ox\cO_X(pD)& \to & 0\\
  &      & \ua &                   &\ua&                    & \ua &      &\\
0 & \to & F_r\ox \cO_X(-qD)& {\to}& V'_{r+s} & {\to} & G_s\ox\cO_X(pD) & \to & 0\\
 &      & \ua &                   &\ua&             & \ua &      &\\
 &  &0 & & 0   &  &0 & \\ \end{array} \]
with $F_r\ox \cO_X(-qD)=i^{-1} V'_{r+s}$ and  $G_s\ox \cO_X(pD)=j(V'_{r+s})$ of ranks $0\leq r\leq p$ and $0\leq s\leq q$ for a subsheaf $V'_{r+s}$ of $V$.

In the following we will discuss the required steps for proving stability of $V$. In total we have to consider all subsheaves $V'_{r+s}$ of $V$ with $0\leq r\leq p$ and $0\leq s\leq q$. However, we can exclude certain cases. 

\underline{\it Step 1:} We first note the cases $(0,0)$, $(p,0)$ and $(p,q)$ do not occur as destabilizing subsheaves. The cases $(0,0)$ and $(p,q)$ are ruled out as we only have to consider subsheaves $V'_{r+s}$ of rank $r+s$ with $0<r+s<p+q$. Note $r=0$ implies $F_r=0$ as $U_p$ does not have a non-zero subsheaf
(the same holds correspondingly for $s=0$). 
Moreover, we can assume \cite[Lemma 4.5]{Friedman book} that the quotient $V/V'_{r+s}$ is torsion free and so cases with $r=p$ need not be considered since $U_p \otimes \cO_X(-qD)) / (F_p\otimes \cO_X(-qD))$ is a torsion sheaf. 

\underline{\it Step 2:} To prove stability of $V$ we need to show that $\mu(V'_{r+s})<0$ for all $0<r+s<p+q$ with $0\leq r \leq p-1$ and $0\leq s \leq q$. We first note
\beqa
(r+s)\mu(V'_{r+s})=(ps-qr)DJ^2 + r\mu(F_r) +s\mu(G_s)\nonumber
\eeqa
and the discussion depends on the sign of $(ps-qr)$. 

\underline{$(ps-qr)<0\colon$} As $DJ^2>0$ by assumption and $\mu(F_r)<0$, respectively, $\mu(G_s)<0$ we get in this case $\mu(V_{r+s}')<0$.

\underline{$(ps-qr)>0\colon$} We have $\mu(F_r)<0$ and $\mu(G_s)<0$ for $0<r<p$ and $0<s<q$.
Further if $s=q$ then $\mu(G_q)\leq 0$ according to \cite[Lemma 4.3]{Friedman book} 
So we get the following subcases:
\begin{enumerate}
\item $\mu(F_r)<0$, \ \ $\mu(G_s)<0$
\item $\mu(F_r)<0$, \ \ $\mu(G_q)<0$
\item $\mu(F_r)<0$, \ \ $\mu(G_q)=0$
\item $r=0$, \ \ $\mu(G_q)<0$
\item $r=0$, \ \ $\mu(G_q)=0$
\end{enumerate}

(i)-(iv)$\colon$ We have to solve $\mu(V_{r+s}')<0$ for $z$, together with $DJ^2>0$ (cf. Proposition \ref{thm22}, Proposition \ref{thm23}, Proposition \ref{thm24}). 

(v)$\colon$ This case has to be excluded as $\mu(V_{0+q}')=(p)DJ^2>0$. 

The following result gives a condition when subsheaves $V'_{0+q}\cong G_q\ox \cO_X(pD)$ do not destabilize $V$ because they do not exist. So the cases
(iv) and (v) would then be excluded if one could show the corresponding assertion about the $f$-map, however, following this line of argumentation we will actually exclude case (v) below (cf. Lemma \ref{fmap}). 

\begin{lem}
\label{21}
Let $U:=U_p\ox \cO_X(-qD)$, $W:=W_q\ox \cO_X(pD)$ and $G:=G_q\ox \cO_X(pD)$.
A sufficient condition for $V$ not to be destabilized by a subsheaf $G$ of $W$ is given by
injectivity of the map
\beqa 
{\Ext}^1(W, U)\stackrel{f}{\to} {\Ext}^1(G,U)\nonumber
\eeqa
\end{lem}

\begin{proof}
We first ask when is it possible that a map $G\to W$ lifts to a map
$G\to V$. To see this consider 
\beqa
\to {\rm Hom}(G, V)\to {\rm Hom}(G, W)\to {\Ext}^1(G,U)\nonumber
\eeqa
showing that the obstruction to lifting an element of ${\rm Hom}(G, W)$
to an element of ${\rm Hom}(G, V)$ lies in ${\Ext}^1(G, U)$.
We have a commutative diagram \[ \begin{array}{ccc}
{\rm Hom}(W, W) &\stackrel{\partial}\to& {\Ext}^1(W, U) \\
\da & & \da \\
{\rm Hom}(G, W)& \to& {\Ext}^1(G, U)
 \end{array} \]
with $\partial(1)=\xi$ the extensions class. So
we conclude a non-zero element of ${\rm Hom}(G, W)$ can be lifted to an element of ${\rm Hom}(G, V)$ exactly when the extension class $\xi$ is in the kernel of \beqa f\colon {\Ext}^1(W, U)\to {\Ext}^1(G, U)\nonumber
\eeqa
thus if $f$ is injective $f(\xi)\neq 0$ and such a lifting does not exist. 
\end{proof}

\underline{$(ps-qr)=0\colon$} This case has in principle to be treated separately in a manner similar to the case $ps-qr>0$.

Below (Section 4, 5, 6) we will show for extensions of type $(p,q)=(n,1)$
the following: the case $ps-qr=0$ does not occur,
the $f$-map arguments for the case $(0,q)$ can be carried through and
the non-split condition can be fulfilled.

Finally, a direct computation gives the Chern classes of $V$
\begin{eqnarray*}
c_1(V)&=&0\\
c_2(V)&=&-\frac{1}{2}pq(p+q)D^2+c_2(U_p)+c_2(W_q)\\
c_3(V)&=&\frac{1}{3}{pq}(p^2-q^2)D^3+2\big( qc_2(U_p)-pc_2(W_q)\big)D+c_3(U_p)+c_3(W_q).
\end{eqnarray*}


\section{Non-Split conditions}

We will now restrict the general set-up of Section 2 to extensions with $U_p$ given by stable pull-back bundles $\pi^*E$ (with $E$ a stable rank $n$ vector bundle on $B$) and $W_q=\cO_X$ such that the resulting vector bundle $V$ has rank $m=n+1$. To prove stability of $\pi^*E$ and $V$ will then be our main focus in the subsequent two sections.

The following result provides sufficient conditions for choosing non-split extensions if $U_p$ is given by $\pi^*E$. 

Let $y:=mx$ and $E_1:=R^1\pi_*\cO_X(-y\si)\ox E\ox \cO_B(-m\al)$ and $E_2:=\pi_*\cO_X(-y\si) \ox \cO_B(-m\al)\ox E$ (the expressions $E_1$ and $E_2$ occur below in the Leray spectral sequence). 

\begin{lem} Let $E$ be a $H$-stable rank $n$ vector bundle on a rational surface $B$ with $c_1(E)=0$, and let $D=x\si+\pi^*\al$ then an extension of $\cO_X(nD)$ by $\pi^*E\ox \cO_X(-D)$ can be chosen non-split for
\begin{enumerate}
\item $x>0\colon$  if $(2H-zc_1)\al \leq 0$ and $\chi(B, E_1)>0$.
\item $x< 0\colon$ if $\chi(B, E_2)<0$.
\item $x=0\colon$  if $\chi(B, E\ox \cO_B(-m\al))<0$.
\end{enumerate}
\label{nspE}
\end{lem}

\begin{proof} We apply the Leray spectral sequence to $\pi\colon X\to B$ and use the projection formula giving 
\beqa
0\to H^1(B, E_2)\to H^1(X, \pi^*E\ox \cO_X(-mD))\to H^0(B, E_1)\to H^2(B,E_2).\nonumber
\eeqa
If $x>0$ then $\pi_*\cO_X(-y\si)=0$ thus
\beqa
H^1(X, \pi^*E\ox \cO_X(-mD))\cong H^0( B, E_1)\nonumber
\eeqa
now Serre duality on $B$ and $[R^1\pi_*\cO_X(-y\si)]^*=\pi_*\cO_X(y\si)\ox K_B^{-1}$ give
\beqa
H^2(B, E_1)^*=H^0(B, \pi_*\cO_X(y\si)\ox \cO_B\big(m\al\big)\ox E^*)\nonumber
\eeqa
where $\pi_*\cO_X(y\si)=\cO_B\oplus K_B^{2}\oplus \cdots \oplus K_B^{y}$ for $y>1$ \cite[Lemma 4.1]{FMWIII} thus
\beqa
H^2(B, E_1)^*=H^0(B,\cO_B\big(m\al\big)\ox E^*)\oplus H^0(B,\bigoplus_{i=2}^yK_B^i\ox\cO_B\big(m\al\big)\ox E^*).\nonumber
\eeqa
Now the first term vanishes if 
\beqa
n\mu\big(\cO_B\big(m\al\big)\ox E^*\big)= (2H-zc_1)\al\leq 0\nonumber
\eeqa
all other terms vanish if 
\beqa
n\mu\big(K_B^i\ox\cO_B\big(m\al\big)\ox E^*\big)=-i(2H-zc_1)c_1+m(2H-zc_1)\al\leq 0\nonumber
\eeqa
as $2H-zc_1\in \C_B$ it follows $(2H-zc_1)c_1\geq 0$ and so we only have to impose
$(2H-zc_1)\al\leq 0$. Note for $B$ the Enriques surface this condition becomes $\al H\leq 0$. For a surface with $K_B^{-1}$ ample and $H=hc_1$, cf. below, the condition becomes $(2h-z)\al c_1\leq 0$; as $2h-z>0$ this is equivalent to $\al c_1\leq 0$, i.e., again $\al H\leq 0$.

Having $H^2(B, E_1)^*=0$, we can now apply the Hirzebruch-Riemann-Roch formula and conclude if $\chi(B, E_1)>0$ then $H^1(X, \pi^*E\ox \cO_X(-mD))$ is non-zero completing the proof of (i).

(ii) and (iii)$\colon$ If $x<0$ then $R^1\pi_*\cO_X(-y\si)=0$ and the Leray spectral sequence gives
\beqa
H^1(X, \pi^*E\ox \cO_X(-mD))\cong H^1(B, E_2)\nonumber
\eeqa
thus if $\chi(B, E_2)<0$ we have $H^1(B, E_2)\neq 0$. 

If $x=0$ then the Leray spectral sequence simplifies (with $\pi_*\cO_X=\cO_B$) \beqa
0\to H^1(B, E\ox \cO_B(-m\al))\to H^1(X, \pi^*E\ox \cO_X(-mD))\to\nonumber
\eeqa
and a sufficient condition for the first space to be non-zero is $\chi(B, E\ox \cO_B(-m\al))<0$.
\end{proof}

Let us state the explicit expressions for $\chi(B,E_1)$, $\chi(B,E_2)$ and $\chi(B, E\ox\cO_B(-m\al))$.
For $y=mx>0$ we note $R^1\pi_*\cO_X(-y\si)=K_B^{1}\oplus K_B^{-1}\oplus\dots \oplus K_B^{1-y}$ for $y>1$ \cite[Lemma 5.16]{FMWIII} and  
\beqa
ch(R^1\pi_*\cO_X(-y\si))=y+A_1c_1+A_2\frac{c_1^2}{2}\nonumber
\eeqa
where we have set $A_1=-1+\frac{y(y-1)}{2}$ and $A_2=1+\frac{y(y-1)(2y-1)}{6}$.
The Hirzebruch-Riemann-Roch formula gives
\beqa
\chi(B, E_1)=y \big(n-c_2(E)+\frac{nm^2}{2}\al^2\big)+A_3\frac{n}{2}c_1^2-A_4nm\al c_1\nonumber
\eeqa
where $A_3=\frac{y(y^2-1)}{3}$ and $A_4=-1+\frac{y^2}{2}$.
For $y=mx<0$ we obtain
\beqa
ch(\pi_*\cO_X(-y\si))=-y-A_1c_1-A_2\frac{c_1^2}{2}\nonumber
\eeqa
and the Hirzebuch-Riemann-Roch formula yields
\beqa
\chi(B, E_2)=-y\Big(n-c_2(E)+\frac{nm^2}{2}\al^2\Big)-A_3\frac{n}{2}c_1^2
+A_4nm\al c_1.\nonumber
\eeqa
If $x=0$ we find
\beqa
\chi(B, E\ox \cO_B(-m\al))=n-c_2(E)+\frac{nm}{2}\al\big(m\al-c_1\big).\nonumber
\eeqa

\begin{rem} Note that for $x>0$ the case $\al = c_1$
does not lead to $H^1(X,\pi^*E\ox\cO_X(-mD))\neq 0$. $E$ being supposed to be stable, one has
$H^0(B, E)=0$; the same holds for the slope zero
stable bundle $K_B\ox E$. However, if $E$ is a $H$-semistable vector bundle of zero slope on the Enriques surface and $H^0(B, E)\neq 0$ then for $x>0$ and $\al=c_1$ one can choose an extension of  $\cO_X(nD)$ by $\pi^*E\ox \cO_X(-D)$ to be non-split.
\end{rem}

\begin{lem}
For $x>0$, the condition $\mu(K_B^{1-y}\ox\cO_B(-m\al)\ox E)>0$ is necessary for the existence of a non-split extension of  $\cO_X(nD)$ by $\pi^*E\ox \cO_X(-D)$. 
\end{lem}
\begin{proof}
From Lemma \ref{nspE} we have for $x>0$ that
\beqa
H^1(X, \pi^*E\ox \cO_X(-mD))\cong H^0( B, R^1\pi_*\cO_X(-y\si)\ox E\ox \cO_B(-m\al)).\nonumber
\eeqa
If  $\mu(K_B^{1-y}\ox\cO_B(-m\al)\ox E)\leq 0$ then we one has $H^0( B, E_1)=0$ and so the extension splits.\end{proof}

\begin{cor}
\label{cor34}
For $x\neq 0$ and $B$ the Enriques surface the following relation necessarily holds when $V$ is stable
$$x\cdot (\al H)<0$$
\end{cor}

\begin{proof}
We apply Lemma \ref{obv}. For $x<0$ from the condition
\beqa
DJ^2=x\si(H-zc_1)^2+z(2H-zc_1)\al\si>0\nonumber
\eeqa
it follows that $H\al>0$ (for the cases with $K_B^{-1}$ ample and $H=hc_1$ we get $\al c_1>0$). 
For $x>0$ the non-split condition gives 
\beqa
n\mu\big(K_B^{1-y}\ox\cO_B\big(-m\al\big)\ox E\big)=(y-1)(2H-zc_1)c_1-m(2H-zc_1)\al>0.\nonumber
\eeqa
becoming $-2m\al H>0$ for the Enriques surface. 
\end{proof}


\section{Stable Extensions on the Enriques CY space}

Let $\pi\colon X\to B$ be an elliptically fibered Calabi-Yau threefold with a section $\si$
and base $B$ given by an Enriques surface, i.e., $h^{1,0}(B)=0$ and $K_B^2={\cO}_B$. 
We first recall some basic properties of these spaces.
$B$ has non-trivial Hodge numbers $h^{1,1}(B)=10, h^{0,0}(B)=h^{2,2}(B)=1$ , so $c_1^2=0$
and $c_2=12$. Further $\phi c_1=0$ for all $\phi \in H^2(B, {\mathbb Z})$
and the intersection form is even \cite{CoDo}. A smooth curve $C$ has $e(C)=-C^2$,
and a generic ('unnodal') $B$ has no smooth rational curves. One gets for the middle cohomology
$H^2(B,{\mathbb Z})={\mathbb Z}^{10}\oplus{\mathbb Z}_2$ with intersection lattice
$$\Gamma^{1,1}\oplus \, E_8^{(-)}={\left( \begin{array}{cc}
0 & 1 \\ 1 & 0 \end{array} \right)}\oplus \, E_8^{(-)}$$
(orthogonal decompositions).
$B$ is always elliptically fibered over $b={\mathbb P^1}$.
However, two of the fibers, $f_1$ and $f_2$, are double fibers:
$f=2f_i$, which prevents $B$ from having a section and
$c_1= f_1 - f_2$ is not effective.

Let us consider the effective cone. On an unnodal $B$
all irreducible curves $C$ have $C^2\geq 0$.
The integral classes in one of the two components
of the cone in $H^2(B, {\mathbb R})$
defined by $C^2\geq 0$ constitute the effective cone (potentially adding
the torsion class $c_1$ does not matter for this if $C\neq 0$; we will not
always mention explicitly this exceptional case).
For $C$ nef (i.e. $DC\geq 0$ for all curves $C$ on $B$)
$|C|$ is base-point-free, and $C$ is ample
if also $C^2\geq 6$ \cite{CoDo}. A $C=xa+yf:=(x,y) \in \Gamma^{1,1}$ is nef
precisely if it is effective in the ${F_0}$-sense, i.e., for $x,y\geq 0$.

Furthermore, we note that $B$ can be represented as the qoutient of a $K3$ surface
by a free involution. The $K3$ can be represented
as a double cover $w^2=f_{4,4}(z_1, z_2)$ of ${\mathbb P^1_{z_1}}\x {\mathbb P^1_{z_2}}$,
branched along a curve of bidegree $(4,4)$,
so elliptically fibered $p_i: K3\ra {\mathbb P^1_{z_i}}$.
The involution is $(z_1, z_2; w)\ra (-z_1, -z_2; -w)$.
This shows also two elliptic fibrations of $B$
with the double fibers over $0$ and $\infty$
\begin{eqnarray*}
K3 & \lra  & B\\
p_{K3} \da & & \; \da \, p\\
{\mathbb P^1_{z_1}}  & \stackrel{(\cdot)^2}{\lra} & {\mathbb P^1_{z_1}} 
\end{eqnarray*}
Note also that in an orbifold limit $T^4/{\mathbb Z_2}$
of $K3$ the involution is $(-1, 1)$ on the complex coordinates $(t_1, t_2)$
of $T^4$, combined with a shift by a half lattice vector in both directions.

The corresponding $\pi_1(B)={\mathbb Z}_2$ is inherited by
the elliptic Calabi-Yau space $X$ which itself is
a quotient by a free involution on $K3\x T^2$ (it acts as described
on $K3$ and as $z\ra -z$ on the $T^2$).
The holomorphic two-form $\Omega_2$ of $K3$ being odd, the holomorphic
three-form $\Omega_2\wedge dz$ is preserved, the quotient $X$
being a Calabi-Yau space of vanishing Euler number. Finally, the second Chern class of $X$ can be obtained by a standard computation (cf. \cite{FMW}) and is given in general by (with $c_i:=c_i(B)$)
\beqa
c_2(X)=\pi^*c_2+11\pi^*c_1^2+12\si \pi^*c_1\nonumber.
\eeqa

Now given a stable vector bundle $E$ of rank $n\geq 2$ with $c_1(E)=0$ on an Enriques surface, we will construct rank $n+1$ vector bundles of trivial determinant on $X$ as non-split extensions of $\cO_X(nD)$ by $\pi^*E_n\ox \cO_X(-D)$ with $D=x\sigma+\pi^*\al$ and prove that $V$ is stable in a region of the K\"ahler cone of $X$.
For this we show first that $\pi^*E$ is stable on $X$ provided that $E$ is stable on the Enriques surface. For the existence of stable vector bundles on Enriques surfaces 
see \cite{kim}, \cite{friedmv}. The fact that $\pi^*E$ is stable on Calabi-Yau threefolds elliptically fibered over the Enriques surface has been used previously in \cite{thom1}.

\begin{lem}
\label{pE}
$\pi^*E$ is (semi-)stable with respect to $J=z\sigma +\pi^*H$ on $X$ if $E$ is (semi-)stable on $B$ with respect to $H$ and with $c_1(E)=0$.
 \end{lem}

\begin{proof} Let $\F$ be a subsheaf of $\pi^*E$ where we can assume that $\pi^*E/\F$ is torsion free \cite[Lemma 4.5]{Friedman book};
so we have $0\to \F\vert_{\si}\to E$ and $c_1(\F\vert_{\si})H<0$ (for semistability $\leq 0$). Similarly we get $0\to \F\vert_{F}\to \cO^r_F$
thus $deg(\F\vert_{F})\leq 0$ as $\cO^r_F$ is semistable (where $r:=rk\, E$). We conclude that $c_1(\F) = -A\si + \la$ with $A\geq 0$ and $\la H < 0$
and $c_1(\F)J^2 = -AH^2\si + 2z \la H \si < 0$ (with $<$ replaced by $\leq$ for semistability).
\end{proof}

In the following Proposition we construct stable bundles $V$ of vanishing first Chern class (recall that $xa<0$ for $x\neq 0$ by Corollary \ref{cor34}).
\begin{prop}
\label{thm22}
Let $V$ be a rank $n+1$ vector bundle on $\pi\colon X\to B$ defined by a non-split extension
$$
0\to \pi^*E\ox \cO_X(-D)\to V\to \cO_X(nD)\to 0
$$
with $E$ an rank $n$, $H$-stable bundle with $c_1(E)=0$ on an Enriques surface $B$ and $D=x\sigma +\pi^*\al$ and $a:= \al H$.
Then $V$ is stable with respect to $J=z\sigma+\pi^*H$ for $|x|< |a|$ and

\begin{enumerate}
\item $x>0$ and $\frac{nx}{1-na}\frac{H^2}{2}<z<\frac{nx}{-na}\frac{H^2}{2}$
\item $x<0$ and $\frac{-nx}{na}\frac{H^2}{2}<z<\frac{-nx}{na-1}\frac{H^2}{2}$
\end{enumerate}
\end{prop}

\begin{proof} We have to consider now the following diagram of exact sequences
\[ \begin{array}{ccccccccc}
  &  &0 & & 0   &  &0 & \\
  &      & \ua &                   &\ua&             & \ua &      &\\
0 & \to & P=\bar{P}\ox \cO_X(-D)& \to & V/V'_{r+1}& \to & T=\bar{T}\ox \cO_X(nD) & \to & 0 \\
  &      & \ua &                   &\ua&             & \ua &      &\\
0 & \to & \pi^*E\ox \cO_X(-D)& \stackrel{i}{\to} & V& \stackrel{j}{\to}
 &\cO_X \ox\cO_X(nD)& \to & 0\\
  &      & \ua &                   &\ua&                    & \ua &      &\\
0 & \to  & F_r\ox \cO_X(-D)& {\to}& V'_{r+1} & {\to} & G_1\ox\cO_X(nD) & \to
 & 0 \\
   &      & \ua &                   &\ua&             & \ua &      &\\
  &  &0 & & 0   &  &0 & \\
 \end{array} \]
In view of the discussion in Section 2 we have to prove stability of $V$ for $
0<r+1<n+1$ with $0\leq r\leq n-1$ and show $\mu(V'_{r+1})<0$. We have
\beqa
(r+1)\mu(V'_{r+1})=(n-r)DJ^2+r\mu(F_r)+\mu(G_1)\nonumber
\eeqa
with $\mu(F_r)<0$ by stability of $\pi^*E$ and $\mu(G_1)\leq 0$ by \cite[Lemma 4.3]{Friedman book}. 
We first note that cases with $(n-r)=0$ do not occur as $r\neq n$. We are left with discussing the cases $(n-r)>0$ leading to the following subcases:
\begin{enumerate}
\item $\mu(F_r)<0$, \ \ $\mu(G_1)<0$
\item $\mu(F_r)<0$, \ \ $\mu(G_1)=0$
\item $r=0$, \ \ $\mu(G_1)<0$
\item $r=0$, \ \ $\mu(G_1)=0$
\end{enumerate}
(ii) and (iii)$\colon$ We must solve the following inequalities for $z$ simultaneously (thereby solving (i))
\beqa
(n-r)DJ^2+r\mu(F_r)<0, \ \ \
(n-r)DJ^2+\mu(G_1)<0, \ \ \ DJ^2>0\nonumber
\eeqa
where the last inequality assures $\mu(\pi^*E\ox \cO_X(-D))<0$.  

By stability of $\pi^*E$ we have $c_1(F_r)=-A\si +\pi^*\la$ with $-A\leq 0$ and $\la H <0$ and $c_1(G_1)=-D_2$ with $D_2=B'\si+\pi^*\bar{\beta}$ an effective divisor. We set $-\bar{\beta}=\beta$ such that $c_1(G_1)=-B'\si +\pi^*\beta$. The slopes of $F_r$ and $G_1$ are given by
\beqa
r\mu(F_r)=-AH^2\si+2z\la H\si, \ \ \
\mu(G_1)=-B'H^2\si+2z\beta H\si\nonumber
\eeqa 
We estimate the first two inequalities gives
\beqa
nDJ^2+r\mu(F_r)<0, \ \ \ nDJ^2+\mu(G_1)<0\nonumber
\label{two}
\eeqa
To evaluate these expressions it is clearly enough to pose the following conditions to the slopes
of $F_r$ and $G_1$
\begin{eqnarray*}
\la H&=&-1 \ \ {\rm and}\ \  A=0, \\
\beta H&=&-1\ \ {\rm and} \ \ B'=0, \ \ \ {\rm resp.}\ \ \ \beta H=0\ \  {\rm and}\ \  -B'=-1
\end{eqnarray*}
where we have chosen for $-B'$, respectively, $\beta H$ the worst case, i.e. the biggest values, which are $-1$. In summary, we have to solve for $z$ the following system (with $a:=\al H$ and all intersection products taken in $B$) 
\begin{eqnarray*}
n(xH^2+2z a)-2 z &<&0\\
n(xH^2+2z a)-H^2 &<&0\\
xH^2+2z a&>&0
\end{eqnarray*}
We find for $x>0$ and $a<0$ the following bounds for $z$ 
\begin{eqnarray*}
z&>& \frac{nx}{1-na}\frac{H^2}{2}\\
z&>& \frac{nx-1}{-na}\frac{H^2}{2}\\
z&<&\frac{-x}{a}\frac{H^2}{2}
\end{eqnarray*}
giving for $0<x<-a$ the condition
\beqa
\frac{nx}{1-na}\frac{H^2}{2}<z<\frac{-x}{a}\frac{H^2}{2}\nonumber
\eeqa
If $x=0$ and $a>0$ we get
\begin{eqnarray*}
na&<&1\\
z&<&\frac{1}{na}\frac{H^2}{2}
\end{eqnarray*}
as $na\ge n$ we find that in case (ii) we cannot solve the conditions which would exclude a destabilizing subsheaf of $V$. 
 
Finally, if $x<0$ and $a>0$ we obtain 
\begin{eqnarray*}
z&<&\frac{nx}{1-na}\frac{H^2}{2}\\
z&<&\frac{1-nx}{na}\frac{H^2}{2}\\
z&>&\frac{-x}{a}\frac{H^2}{2}
\end{eqnarray*}
giving for $-a<x<0$ the condition
\beqa
\frac{-x}{a}\frac{H^2}{2}<z<\frac{nx}{1-na}\frac{H^2}{2}\nonumber
\eeqa

(iv)$\colon$ This case will be treated in Lemma \ref{fmap} below. We show that potential subsheaves
of $V$ of type $V'_{0+1}$ with $\mu(V'_{0+1})=nDJ^2>0$ do not exist. 
\end{proof}
To treat the $(0,1)$ cases with $\mu(G_1)=0$ let us first determine the general structure of $G_1$.

\begin{lem}
$G_1$ has the structure $G_1=\cO_X(-D_2) \otimes I_Y$
with $D_2\geq 0$ and codim $Y\geq 2$.\\
\end{lem}
\begin{proof}
In the $(0,1)$ case $G_1\ox \cO_X(nD)=V'\hookrightarrow V_2$ and the torsion sheaf $\bar{T}$ is a quotient of $\cO$, so $\bar{T}=\cO_Z$ for a subscheme $Z$ of $X$, and, considering the vertical
sequence on the right, $G_1=I_Z$ where $I_Z$ is the ideal  sheaf of $Z$. Further $D_2J^2 \ge 0$ for $D_2=c_1(I_Z)$ which comprises the codim $=1$ components of $Z$; so
$I_Z = \cO_X(-D_2)\otimes I_Y $
where $Y$ $j:Y\ra X$ is the closed immersion 
of the union of all components  of $Z$ of codim $\geq 2$. \end{proof}

This result and the following result applies for any extension of type 
\beqa
0\to V_n\ox \cO_X(-D)\to V\to \cO_X(nD)\to 0\nonumber
\eeqa
on a general Calabi-Yau threefold $X$. 
\begin{lem}
\label{fmap}
$V'=I_Y\ox \cO_X(nD)$ with codim $Y\geq 2$ does not occur as a subsheaf of $V$.
\end{lem}
\begin{proof}
Let ${\bar V}_n=V_n\ox \cO_X(-D)$. By Lemma \ref{21} we have to prove that 
$f\colon {\Ext}^1(\cO_X(nD), \bar{V}_n)\to {\Ext}^1(I_Y(nD), \bar{V}_n)$ is injective. 
For this consider 
\beqa
0\to I_Y(nD)\to \cO_X(nD)\to  \cO_Y(nD)\to 0\nonumber
\eeqa
taking ${\rm Hom}$$(\cdot, \bar{V}_n)$ yields
\begin{eqnarray*}
0&\to& {\rm Hom}(\cO_X(nD),\bar{V}_n)\to {\rm Hom}(I_Y(nD),\bar{V}_n)\to {\Ext}^1(\cO_Y(nD),\bar{V}_n)\to\\
&\to&{\Ext}^1(\cO_X(nD),\bar{V}_n)\stackrel{f}\to {\Ext}^1(I_Y(nD),\bar{V}_n)\to \cdots
\end{eqnarray*}
We have to show that ${\Ext}^1(\cO_Y(nD),\bar{V}_n)=0$. Now Serre duality gives 
${\Ext}^1(\cO_Y(nD),{\bar V}_n)\cong {\Ext}^2(\bar{V}_n,\cO_Y(nD))^*$ further we have
\beqa
{\Ext}^2(\bar{V}_n,\cO_Y(nD))^*={\Ext}^2(\cO_X,\bar{V}^*_n\ox \cO_Y(nD))^*=H^2(X,\bar{V}_n^*\ox\cO_Y(nD))^*\nonumber
\eeqa
now as $\cO_Y\equiv j_*\cO_Y$ we have $H^2(X,\bar{V}_n^*\ox\cO_Y(nD))^*=
H^2(X,j_*\big(j^*(\bar{V}_n^*\ox\cO_X(nD))\big))^*=H^2(Y,j^*(\bar{V}_n^*\ox\cO_X(nD)))^*=0$
for codim $Y\geq 2$.
\end{proof}

\begin{rem}
As for $x=0$ we cannot assure solvability of the numerical slope conditions,
one would need to give a condition such that
$$0\to F_r\ox O_X(-D)\to V'_{r+1}\to I_Y\ox \cO_X(nD)\to 0$$ for 
$0<r<n$ does not occur as potential subsheaf of $V$.
\end{rem}


\section{Stable extensions on CY spaces with del Pezzo surface base}

In this section we will consider elliptically fibered Calabi-Yau threefolds $\pi\colon X\to B$ whose base has an ample $K_B^{-1}$. As in the previous section we consider extensions $V$ of $\cO_X(nD)$ by $\pi^*E\ox\cO_X(-D)$. We first note

\begin{lem}
$\pi^*E$ is (semi-)stable on $X$ with respect to $J=z\si +\pi^*H\in {\C}_X$ $($i.e. $H-zc_1\in {\C}_B$, so $z<h)$
if $E$ is (semi-)stable on $B$ with respect to $H=hc_1$ and has $c_1(E)=0$. \end{lem}

\begin{proof}
Let $\F$ be a subsheaf of $\pi^*E$ where we can assume that $\pi^*E/\F$ is torsion free \cite[Lemma 4.5]{Friedman book};
so we have $0\to \F\vert_{\si}\to E$ and $c_1(\F\vert_{\si})H<0$ (for semistability $\leq 0$). Similarly we get $0\to \F\vert_{F}\to \cO^r_F$
thus $deg(\F\vert_{F})\leq 0$ as $\cO^r_F$ is semistable (where $r:=rk\, E$). Then for $H-zc_1\in {\C}_B$ and
$c_1(\F) = -A\si + \la$ with $A\geq 0$ and $\la H \leq (Ac_1+\la)H< 0$
(the latter and the following $<$ are $\leq$ for semistability)
\beqa
\label{general decomposition}
c_1(\F)J^2 = -A(H-zc_1)^2\si + z ( 2H-zc_1)\la\si < 0.\nonumber
\eeqa
\end{proof}
We can now proceed as in the previous section and prove stability of the extension $V$.
\begin{prop}
\label{thm23}
Let $V$ be a rank $n+1$ vector bundle on $\pi\colon X\to B$ defined by a non-split extension
$$
0\to \pi^*E\ox \cO_X(-D)\to V\to \cO_X(nD)\to 0
$$
with $E$ an rank $n$ bundle with $c_1(E)=0$, stable with respect to $H=hc_1$, and $D=x\sigma +\pi^*\al$. Then $V$ is stable with respect to $J=z\sigma+\pi^*H$ for $0<|x|< |a|$ and for 
$z$ in the following ranges $(\zeta:=h-z)$
\begin{enumerate}
\item $0<x<-a$ and $\frac{nx}{n(xc_1^2-a)+1}H^2<h^2-\zeta^2<\frac{nx}{n(xc_1^2-a)}H^2$
\item $-a<x<0$ and $\frac{nx}{n(xc_1^2-a)}H^2<h^2-\zeta^2<\frac{nx}{n(xc_1^2-a)+1}H^2$
\end{enumerate}
\end{prop}

\begin{proof}
The proof is completely parallel to the proof of Proposition \ref{thm22} up to the following consideration.
The slopes of $F_r$, $G_1$ and the expression for $DJ^2$ are given by
\begin{eqnarray*}
r\mu(F_r)&=&-A\big(h-z\big)^2c_1^2\si+z\big(2h-z)\la c_1\si\\
\mu(G_1)&=&-B'\big(h-z\big)^2c_1^2\si+z\big(2h-z)\beta c_1\si\\
DJ^2&=&x(h-z)^2c_1^2\si+z\big(2h-z)\al c_1\si
\end{eqnarray*}
inserting these expressions in the estimated inequalities
\beqa
nDJ^2+r\mu(F_r)<0, \ \ \ nDJ^2+\mu(G_1)<0\nonumber
\eeqa
and imposing as in Proposition \ref{thm22} the following conditions to the slopes
of $F_r$ and $G_1$
\begin{eqnarray*}
\la c_1&=&-1 \ \ {\rm and}\ \  A=0, \\
\beta c_1&=&-1\ \ {\rm and} \ \ B'=0, \ \ \ {\rm resp.}\ \ \ \beta c_1=0\ \  {\rm and}\ \  -B'=-1.
\end{eqnarray*}
we find the following inequalities, which have to be solved for $z$ (here and in the following 
all intersection products are taken in $B$)
\begin{eqnarray*}
nxh^2c_1^2+\big(na-1-nxc_1^2\big)\big(2h-z\big)z&<&0\\
\big(nx-1\big)h^2c_1^2+\big(na+c_1^2-nxc_1^2\big)\big(2h-z\big)z&<&0\\
xh^2c_1^2+\big(a-xc_1^2\big)\big(2h-z\big)z&>&0
\end{eqnarray*}
Define $\zeta:=h-z$ such that $0<\zeta<h$.
For $x>0$ and $a<0$ we get
\begin{eqnarray*}
h^2-\zeta^2&>&\frac{nxh^2c_1^2}{n(xc_1^2-a)+1}\\
h^2-\zeta^2&>&\frac{(nx-1)h^2c_1^2}{n(xc_1^2-a)-c_1^2}\\
h^2-\zeta^2&<&\frac{xh^2c_1^2}{xc_1^2-a}
\end{eqnarray*}
so we get for $0<x<-a$ 
\beqa
\frac{nxh^2c_1^2}{n(xc_1^2-a)+1}<h^2-\zeta^2<\frac{xh^2c_1^2}{xc_1^2-a}.\nonumber
\eeqa
For $x=0$ and $a>0$ we find from the first inequality above
\beqa
(na-1)(h^2-\zeta^2)<0\nonumber
\eeqa
but $(h^2-\zeta^2)>0$ and $(na-1)>0$.

Finally, for $x<0$ and $a>0$ we get
\begin{eqnarray*}
h^2-\zeta^2&<&\frac{nxh^2c_1^2}{n(xc_1^2-a)+1}\\
h^2-\zeta^2&<&\frac{(nx-1)h^2c_1^2}{n(xc_1^2-a)-c_1^2}\\
h^2-\zeta^2&>&\frac{xh^2c_1^2}{xc_1^2-a}
\end{eqnarray*}
for $-a<x<0$ we get
\beqa
\frac{xh^2c_1^2}{xc_1^2-a}<h^2-\zeta^2<\frac{nxh^2c_1^2}{n(xc_1^2-a)+1}.\nonumber
\eeqa
\end{proof}


\section{Stable extensions of spectral cover bundles}

In this section we will study non-split extensions of stable spectral cover bundles $V_n$
on $\pi\colon X\to B$ with $B$ either given by a Hirzebruch surface (or blow-ups of it), a del Pezzo surface or an Enriques surface. We first recall the notion stable spectral cover bundle. 

Let $X$ be an elliptically fibered Calabi-Yau threefold with a section $\si$, let $C$ be an irreducible surface in the linear system $|n\si+\eta|$ and $i\colon C\to X$ the immersion of $C$ into $X$ and let $L$ be a rank one sheaf on $C$. We say $V_n$ is a spectral cover bundle of rank $n$ if $V_n=\pi_{1*}(\pi_2^*(i_*L)\ox \mathcal{P})$ where $\mathcal{P}$ is the Poincar\'e sheaf on the fiber product $X\times_B X$ and $\pi_{1,2}$ are the respective projections on the first and second factor.
Moreover, $V_n$ is stable with respect to $J=\epsilon J_0+\pi^*H$ for $0<\epsilon <\epsilon_0$ as stated in Theorem 7.1 in \cite{FMWIII} (we will always assume $\epsilon$ sufficiently small).  Furthermore, note that various aspects of the spectral cover construction have been studied in 
\cite{bj},\cite{C},\cite{DOPWII},\cite{DOPWIII},\cite{AH1}.

Let $H$ be an ample divisor in $B$, we define the minimal $H$-degree as follows
\beqa
(\Lambda H)_{\rm min}={\rm min}\{\Lambda\cdot H\vert\  \Lambda \in H^2(B,\mathbb{Z}) \ {\rm effective},\  \Lambda\cdot H>0\}\nonumber
\eeqa
which will be useful for defining the minimal slope of the subsheaf $G_1$.

\begin{prop}
\label{thm24}
Let $V_n$ be a stable spectral cover bundle of rank $n$ on $\pi\colon X\to B$ and let $V$ be defined by a non-split extension $$
0\to V_n\ox \cO_X(-D)\to V\to \cO_X(nD)\to 0
$$
with $D=x\si+\pi^*\al$ then $V$ is stable with respect to $J=\epsilon \si+\pi^*H$ for $x=0$ and
$0 < n\al {H} < (\Lambda {H})_{\rm min}.$
\end{prop}

\begin{proof}
As in Proposition \ref{thm22} we have to treat the cases (i)-(iv). 

(i)-(iii)$\colon$ By stability of $V_n$ we have $c_1(F_r)=-A\si +\pi^*\la$ with $A> 0$ by Theorem 7.1 of \cite{FMWIII} and $c_1(G_1)=-D_2$ with $D_2=B'\si+\pi^*\bar{\beta}$ an effective divisor. We set $-\bar{\beta}=\beta$ such that $c_1(G_1)=-B'\si +\pi^*\beta$. The slopes of $F_r$ and $G_1$ and the expression for $DJ^2$ are given by (we define $\bar{H}:=2H-\epsilon c_1$)
\begin{eqnarray*}
r\mu(F_r)&=&-AH^2\si+\epsilon\big(\la+Ac_1\big)\bar{H}\si\\
\mu(G_1)&=&-B'H^2\si+\epsilon\big(\beta+Bc_1\big)\bar{H}\si\\
DJ^2&=&xH^2\si+\epsilon\big(\al-xc_1\big)\bar{H}\si.
\end{eqnarray*}
As in Proposition \ref{thm23} we have to consider $nDJ^2+r\mu(F_r)<0$ which becomes
\beqa
\big(nx-A\big)H^2\si+\epsilon\big(n\al+\la-(nx-A)c_1\big)\bar{H}\si<0.\nonumber
\eeqa
For $nDJ^2+\mu(G_1)<0$ we have to consider for $-B'$ and $\beta \bar{H}$ the worst case, i.e., when $|\mu(G_1)|$ is minimal. This will be achieved for either $B'=0$ and $-\beta \bar{H}=(-\beta \bar{H})_{\rm min}>0$ or $-B'=-1$ and $\beta \bar{H}=0$. 
\begin{eqnarray*}
n x H^2\si+\epsilon n\big(\al-xc_1\big)\bar{H}\si-\epsilon(-\beta \bar{H})_{\rm min}\si&<&0\\
\big(nx-1\big)H^2\si+\epsilon\big(n\al-(nx-1)c_1\big)\bar{H}\si&<&0
\end{eqnarray*}

For $x>0$ we cannot solve the first inequality as we have to assume $A=1$ as the worst case. The cases with 
$x<0$ are ruled out by the positivity condition $DJ^2>0$. 

If $x=0$ then $DJ^2=\epsilon(\al \bar{H})>0$ implies $\al \bar{H}>0$ and we get (intersections taken in $B$)
\begin{eqnarray*}
-AH^2+\epsilon \big(n\al+\la+Ac_1\big)\bar{H}&<&0\\
n \al \bar{H}-(-\beta\bar{H})&<&0\\
-H^2+\epsilon(n\al+c_1)\bar{H}&<&0
\end{eqnarray*}
where the second constraint implies $n\al \bar{H}<(-\beta \bar{H})_{\rm min}$.

Finally, the case (iv) is ruled out by Lemma \ref{fmap}.
\end{proof}

Let us give an example for solving $0<n\al \bar{H}<(-\beta \bar{H})_{\rm min}$.
\begin{exa}
Let $B$ be the Enriques surface then $\bar{H}=2H$; further let $H, \al, \beta \in \Gamma^{1,1}$.
We fix a polarization $H=(v,v+1)$ and set $\beta=-\bar{\beta}=-(e,f)$. We take $\al=(1,-1)$ such that
$\al H=1$ and so we have to solve the inequality for $n$. Then $(-\beta H)_{\rm min}=v$ and we get solutions for $0<n<v$. 
\end{exa}

Finally, we state the conditions such that an extension of $\cO_X(n\pi^*\al)$ by $V_n\ox \cO_X(-\pi^*\al)$
can be chosen non-split. 

\begin{lem}
Let $V_n$ be a stable spectral cover bundle of rank $n$ and let $\al\in H^2(B,\mathbb{Z})$ then an extension of $\cO_X(n\pi^*\al)$ by $V_n\ox \cO_X(-\pi^*\al)$ can be chosen non-split if $
\chi(\mathcal{A}, R^1\pi_*V_n\vert_{\mathcal{A}}\ox \cO_B(-m\al)\vert_{\mathcal{A}})>0
$ with $\mathcal{A}:=\pi(C\cap \si)$ and $C\in |n\si+\pi^*\eta|$.
\end{lem}

\begin{proof}
Applying the Leray spectral sequence to $\pi\colon X \to B$ yields
\begin{eqnarray*}
0&\to& H^1(B, \pi_*V_n \ox \cO_B(-m\al))\to H^1(X, V_n\ox \cO_X(-m\pi^*\al))\\
&\to& H^0(B,R^1\pi_*V_n\ox \cO_B(-m\al))\to H^2(B, \pi_*V_n \ox \cO_B(-m\al))\to
\end{eqnarray*}
For a given spectral cover bundle $V_n$
one has $\pi_*V_n=0$. At a generic point $b\in B$ one has the stalk
$(\pi_*V_n)_b=H^0(F, V_n|_F)=\bigoplus_{i=1}^nH^0(F,{\cO}_F(q_i-p))$
where $p=\si F$ is the zero element in the group law on the fibre
$F$ over $b\in B$ and $q_i$ are the points at which the spectral cover
of $V_n$ intersects $F$.
Now ${\cO}(q_i-p)$ is generically a non-trivial bundle of degree zero which
over an elliptic curve admits no global sections.
Thus $H^0(F,{\cO}_F(q_i-p))=0$ for all $i$ and so $(\pi_*V_n)|_b=0$.
However, since $V_n$ is torsion free, $\pi_*V_n$ is also torsion free.
Thus $(\pi_*V_n)|_b=0$ for generic $b\in B$ gives $\pi_*V_n=0$ everywhere. It follows
\beqa
H^1(X, V_n\ox \cO_X(-m\pi^*\al))\cong H^0(B,R^1\pi_*V_n\ox \cO_B(-m\al))\nonumber
\eeqa
The sheaf $R^1\pi_*V_n$ has support on $\mathcal{A}=\pi(C\cap \si)$ of class $\eta-nc_1$ in $B$ and $H^0(B, R^1\pi_*V_n\ox \cO_B(-m\al))\cong H^0(\mathcal{A}, R^1\pi_*V_n\vert_{\mathcal{A}}\ox \cO_B(-m\al)\vert_{\mathcal{A}})$.
The Grothendieck-Riemann-Roch theorem for $\pi\colon X\to B$ gives 
\beqa
c_1(R^1\pi_*V_n)=\eta-nc_1\nonumber
\eeqa
The Riemann-Roch formula gives 
\beqa
\chi(\mathcal{A}, R^1\pi_*V_n\vert_{\mathcal{A}}\ox \cO_B(-m\al)\vert_{\mathcal{A}})=\frac{3}{2}(\eta-nc_1)^2-m\al(\eta-nc_1)\nonumber
\eeqa
so if $\frac{3}{2}(\eta-nc_1)^2-m\al(\eta-nc_1)>0$ then $H^0(S, R^1\pi_*V_n\vert_{\mathcal{A}}\ox \cO_B(-m\al)\vert_{\mathcal{A}})\neq 0$
and the extension can be chosen to be non-split. Note that irreducibility of the spectral surface (required for stability of $V_n$) demands $\eta\geq nc_1$. 
\end{proof}


\section{Physical Solutions}

A compactification of the perturbative $E_8\times E_8$ heterotic string on a Calabi-Yau threefold $X$
requires mathematically to construct a pair of stable holomorphic vector bundles ($V_1,V_2$) of the same slope and trivial determinant on $X$. Consistency of the physical theory requires the bundles to satisfy the topological constraint
\beqa
\label{anom}
c_2(TX)=\sum_{i}c_2(V_i),\nonumber
\eeqa
the necessary condition for solutions to the heterotic anomaly condition. The spectral cover construction does not lead to vector bundles which solve this topological constraint as $[W]:=c_2(TX)-\sum_i c_2(V_i)$ is non-zero. If $[W]$ is non-zero then physically one expects a five-brane to contribute a source term $\delta_4$ (a current that integrates to one in the direction transverse to a single five-brane) to the Bianchi identity for the three-form $H$. To each five-brane one associates such a four-form delta function source. The class $[W]$ is then the Poincar\'e dual of an integer sum of all these sources and thus $[W]$ should be integral, representing a class in $H_2(X,{\mathbb Z})$. We will use the same expression $[W]$ for an integral homology class in $H_2(X,{\mathbb Z})$, an integral cohomology class in $H^4(X,{\mathbb Z})$ and the de Rham cohomology  
class in $H^4_{\rm DR}(X,{\mathbb R})$ (i.e., as $H^p(X,{\mathbb Z})\to H^p_{\rm DR}(X,{\mathbb R})$ is not injective, integral classes are identified with the images of $H^p(X,{\mathbb Z})$ in $H^p(X,{\mathbb R})$). $[W]$ can be further specified taking into account
that supersymmetry requires that five-branes are wrapped on holomorphic curves thus $[W]$ 
must correspond to the homology class of holomorphic curves. Algebraic classes include negative classes, however, these lead to negative charges, which are unphysical, and so they have to be excluded constraining $[W]$ to be an effective class. 
Thus for a given Calabi-Yau threefold $X$ the effectivity of $[W]$ constrains the choice of
vector bundles $V$. Consideration of the physical background of a heterotic string compactification
\cite{ACu06} reveals that $[W]=0$ is the case favorite by consistency requirements. 

What we will show now is that the vector bundles constructed above allow to solve the topological 
constraint with $[W]=0$. To solve the anomaly constraint we set $V_1=V$ and $V_2=0$ we find
\beqa
c_2(X)-c_2(V)=\pi^*w_B\, \si+a_f [F]\nonumber
\eeqa
where $[F]$ denotes the class of a fiber of $X$. Let $V$ be given by 
\beqa
0\to \bar{V}_n\otimes \cO_X(-D)\to V\to \cO_X(nD)\to 0\nonumber
\eeqa
where $\bar{V}_n$ will be specified below as either $\pi^*E$ or a spectral bundle $V_n$, we find
(set $c_2(\bar{V}_n)=\pi^*\phi\, \si+\pi^*\omega$ with $\phi\in H^2(B, \mathbb{Z})$ and $\omega\in H^4(B,\mathbb{Z})$)
\begin{eqnarray*}
w_B&=&12c_1-\phi+\frac{1}{2}n(n+1)x(2\al-xc_1)\\
a_f&=&c_2+\frac{1}{2}n(n+1)\al^2 - \omega+11c_1^2.
\end{eqnarray*}

\subsection{The case of $B$ the Enriques surface} Let $\bar{V}_n=\pi^*E$, if $B$ is an Enriques surface the problem of finding solutions to $[W]=0$ simplifies as the following 
result shows.
\begin{prop}
Let $\pi\colon X\to B$ with $B$ an Enriques surface then the physical constraint $w_B\geq 0$ implies $x=0$.
\end{prop}

\begin{proof}
For Enriques base we find
\begin{eqnarray*}
\label{w_B condition}
w_B&=&\frac{1}{2}n(n+1)x(2\al-xc_1)\geq 0\\
a_f&=&12+\frac{1}{2}n(n+1)\al^2-c_2(E) \geq 0
\end{eqnarray*}
to solve $w_B\geq 0$ requires $x\al\geq 0$ (we already argued from the non-split
condition that $\al \neq 0, c_1\;$). So $x\al H \geq 0$ contradicting Lemma \ref{nspE}.
\end{proof}

If $V_2$ is non-trivial the argument remains valid as $w_B=\sum_{i=1}^2 a_i\, x_i\al_i\geq 0$
(with $a_i>0, a_1=\frac{n(n+1)}{2}$) gives $w_B H\geq 0$ in contradiction to Lemma \ref{nspE}.

But $x=0$ is the case where the existence of stable bundles could not be assured above. 

\subsection{The case of $B$ a del Pezzo surface}

Let $\pi\colon X\to B$ an elliptic Calabi-Yau threefold with $K_B^{-1}$ ample. Let $\bar{V}_n=\pi^*E$. 
In contrast to the case of the Enriques base it is now possible to satisfy $w_B\geq 0$ while having $x\neq 0$. One finds $[W]=0$ for the choices
\begin{eqnarray*}
\al=\Big(\frac{x^2}{2}-\frac{12}{n(n+1)}\Big)\frac{c_1}{x}
& \Longrightarrow  & w_B=0\\
c_2(E)=c_2+11c_1^2+\frac{n(n+1)}{2}\al^2
& \Longrightarrow  & a_f=0
\end{eqnarray*}

For $x>0$ the non-split condition (cf. Lemma \ref{nspE}, (i)) is satisfied if (set $m=n+1$)
\beqa 
x^2\leq\frac{24}{nm}\nonumber
\eeqa
\beqa
2n+\Big(\frac{(3m^3+m^2)nx^2}{12}+\frac{144}{mx}-\frac{37}{3}n-20\Big)c_1^2>24.\nonumber
\eeqa
For instance, for building an $SO(10)$ GUT model without fivebranes one can use
the twist $D=\si -\pi^*c_1/2$ and a rank $n=3$ bundle $E$
on a base ${F_r}$ of $c_2(E)=104$. For another case one may 
construct an $E_6$ GUT model without fivebranes
from using the twist $D=2\si$ and a plane bundle of $c_2(E)=92$.

\subsection{Extensions by spectral bundles}
Let $\pi\colon X\to B$ an elliptic Calabi-Yau threefold with $B$ either given by a Hirzebruch surface (or blow-ups of it), a del Pezzo surface or an Enriques surface. Let 
$\bar{V}_n=V_n$ a spectral rank $n$ vector bundle has second Chern class equals to \cite{FMW}
\beqa
c_2(V_n)=\pi^*\eta\, \si-\frac{1}{24}\pi^*c_1^2(n^3-n)+\frac{1}{2}(\la^2-\frac{1}{4})n\pi^*\eta(\pi^*\eta-n\pi^*c_1)\nonumber
\eeqa
The condition $c_1(V_n)=0$ imposes constraints on the spectral data \cite{FMW}. One finds: if $n$ is even then $\la=m+\frac{1}{2}$ and $m\in\mathbb{Z}$. If $n$ is odd one has
$\la=m$ and $\eta\equiv c_1\ ({\rm mod}\ 2)$. For this set-up we get
\begin{eqnarray*}
w_B&=&(12c_1-\eta)\\
a_f&=&c_2+11c_1^2+\frac{1}{2}n(n+1)\al^2+\frac{1}{24}(n^3-n)c_1^2-\frac{1}{2}(\la^2-\frac{1}{4})n\eta(\eta-nc_1)
\end{eqnarray*}
Now $w_B=0$ is solved for $\eta=12c_1$. Then for $a_f=0$, we have to solve
\beqa
c_2+c_1^2\Big(11+\frac{n^3-n}{24}-\frac{1}{2}(\la^2-\frac{1}{4})(12-n)n\Big)+\frac{1}{2}n(n+1)\al^2=0\nonumber
\eeqa
To give an example let us assume $B=F_0$. We get for instance $a_f=0$ for $n=2$ and $m=1$ and $\al=(1,-11)$. The bundle is stable for $H=(3,34)$ so $\al H=1$ and $0<n \al H<(-\beta H)_{\rm min}$ becomes $0<2<3$ is satisfied. The non-split condition $\frac{3}{2}(\eta-nc_1)^2-m\al(\eta-nc_1)>0$ is  satisfied.


\begin{thebibliography}{mm}

\bibitem{FMW}
R. Friedman, J.W. Morgan and E. Witten,
{\em Vector Bundles And F Theory},
hep-th/9701162, Commun.Math.Phys. {\bf 187} (1997) 679.


\bibitem{FMWIII}
R. Friedman, J.W. Morgan and E. Witten,
{\em Vector Bundles over Elliptic Fibrations},
Jour. Alg. Geom. {\bf 8} (1999) 279-401, alg-geom/9709029.

\bibitem{Don1}
R. Donagi, {\em Principal bundles on elliptic fibrations}, Asian J. Math. {\bf 1} (1997), 214-223, alg-geom/9702002.

\bibitem{ACu06}
B. Andreas and G. Curio, {\em Heterotic models without fivebranes}, to appear. 

\bibitem{Bfibre}
B. Andreas and G. Curio,
{\em Standard Models from Heterotic String Theory},
hep-th/0602247.

\bibitem{thom1}
R. P. Thomas, {\em Examples of bundles on Calabi-Yau 3-folds for string theory compactifications}, Adv. Theor. Math. Phys. {\bf 4}, 231-247, 2000.

\bibitem{CoDo}
F. Cossec and I. Dolgachev, {\em Enriques Surfaces I},
Birkh\"auser 1989.

\bibitem{kim}
H. Kim, 
{\em Moduli Spaces of Stable Vector Bundles on Enriques Surfaces},
Nagoya Math. J., Vol. {\bf 150} (1998) 85.

\bibitem{friedmv}
R. Friedman, 
{\em Rank two vector bundles over regular elliptic surfaces}, Invent Math., 96 (1989), 283-332.

\bibitem{Friedman book}
R. Friedman, {\em Algebraic Surfaces and Holomorphic Vector Bundles},
Springer, Universitext (1998).

\bibitem{bj}
B. Andreas,
{\em On Vector Bundles and Chiral Matter in N=1 Heterotic Compactifications},
hep-th/9802202, JHEP {\bf 9901} (1999) 011.

\bibitem{C}
G. Curio, {\em Chiral matter and transitions in heterotic string models},
hep-th/9803224, Phys.Lett. {\bf B435} (1998) 39.


\bibitem{DOPWII}
R. Donagi, B. Ovrut, T. Pantev and D. Waldram,
{\em Standard-model bundles},
math.AG/0008010, Adv.Theor.Math.Phys. {\bf 5} (2002) 563.

\bibitem{DOPWIII}
R. Donagi, B. Ovrut, T. Pantev and D. Waldram,
{\em Spectral involutions on rational elliptic surfaces},
math.AG/0008011, Adv.Theor.Math.Phys. {\bf 5} (2002) 499.

\bibitem{AH1}
 B.. Andreas and  D. Hern\'andez Ruip\'erez, 
 {\em U(n) Vector Bundles on Calabi-Yau Threefolds
for String Theory Compactifications}, hep-th/0410170     Adv. Theor. Math. Phys. {\bf 9}, 253, 2006.\\

\end{thebibliography}
\end{document}